\documentclass[12pt, a4paper]{amsart}
\usepackage{amssymb}
\usepackage{amsmath}
\usepackage[english]{babel}
\usepackage{graphicx}
\newtheorem{theorem}{Theorem}[section]
\newtheorem{lemma}[theorem]{Lemma}

\newtheorem{aux}[theorem]{Auxiliary lemma}

\theoremstyle{definition}
\newtheorem{definition}[theorem]{Definition}
\newtheorem{example}[theorem]{Example}
\newtheorem{question}[theorem]{Question}

\newtheorem{remark}[theorem]{Remark}
\newtheorem{remarks}[theorem]{Remarks}

\newtheorem{examples}[theorem]{Examples}
\theoremstyle{remark}
 \numberwithin{equation}{section}

 \DeclareMathOperator{\re}{Re}
\DeclareMathOperator{\im}{Im}

\renewcommand{\phi}{\varphi}

\begin{document}


\title{More approximation on disks}

\author{Peter de Paepe}
\author{Jan Wiegerinck}
\address{Korteweg -- de Vries Institute, University of Amsterdam,
Plantage Muidergracht 24, 1018 TV, Amsterdam, The Netherlands}
\email{depaepe@science.uva.nl} \email{janwieg@science.uva.nl}
\date{December 5, 2005}

\begin{abstract} In this paper we study the function algebra generated by $z^2$ and $g^2$ on a small
closed disk centered at the origin of the complex plane. We prove,
using a biholomorphic change of coordinates and already developed
techniques in this area, that for a large class of functions $g$ this algebra consists of
all continuous functions on the disk.
\end{abstract}

\maketitle

\setlength{\parindent}{0pt} \addtolength{\textheight}{45pt}
\DeclareGraphicsRule{*}{eps}{*}{}
 \setcounter{page}{1}
 \section{Introduction}
 Let  $g$ be a $C^1$ function defined in a
neighbourhood of the origin in the complex plane, 
with $g(0)=0$, $g_z(0)=0$, $g_{\bar z}(0)=1$  (i.e. $g$ behaves like
$\bar z$ near $0$), and such that $z^2$ and $g^2$ separate points
near $0$. Is it possible to find a small closed disk 
 $D$ about 0 in the complex plane,
so that every
continuous function on $D$ can be approximated uniformly on $D$ by polynomials in
$z^2$ and $g^2$? In other words is the function algebra $[z^2,
g^2;D]$ on $D$ generated by $z^2$ and $g^2$, i.e. the uniform
closure in $C(D)$ of the polynomials in $z^2$ and $g^2$, equal to
$C(D)$?  It has been shown that both answers {\it no} and  {\it yes } 
are possible, cf. \cite{[8],[5]}.
\\
The motivating question for this approximation problem was whether
$[z^2, \bar z^2 +\bar z^3;D]$ equals $C(D)$. The answer has been
given recently by O'Farrell and Sanabria-Garc\'ia  and  is $no$,
cf. \cite{[6]}.
\\
 The crucial point in showing whether or not  the algebra $[z^2, g^2 ; D]$  coincides with $C(D)$, is to determine
 whether or not the
preimage of $X=(z^2, g^2)(D)$
 under the map $\Pi (\zeta _1, \zeta _2 )=(\zeta
_1^2, \zeta _2^2 )$ is polynomially convex. Now the set
${\Pi}^{-1}(X)$ consists of the following four disks:
\begin{align*}
D_1&=\{(z, g(z)) \, : \, z\in D\},
\\
D_2&=\{(-z, -g(z)) \, : \, z\in D\} = \{(z, -g(-z)) \, : \, z\in
D\},\\
D_3&=\{(-z,g(z)) \, : \, z\in D\}, \\
D_4&=\{(z, -g(z)) \, : \, z\in D\}=\{(-z, -g(-z)) \, : \, z\in
D\}.
\end{align*}
 In this situation our problem boils down to the (non-)polynomial
convexity of
 $D_1\cup D_2$.
\\
\\
 An appropriate tool in this context is Kallin's lemma: suppose $X_1$ and $X_2$
are polynomially convex subsets of $\mathbb {C}^n$, suppose there is
a polynomial $p$ mapping $X_1$ and $X_2$ into two polynomially
convex subsets $Y_1$ and $Y_2$ of the complex plane such that 0 is
a boundary point of both $Y_1$ and $Y_2$ and with $Y_1\cap
Y_2=\{0\}$. If $ \; p^{-1}(0)\cap (X_1\cup X_2)$ is  polynomially
convex,
 then $X_1\cup X_2$ is polynomially
convex, \cite{[1],[10]}.
\\
In \cite{[3]} Nguyen and the first author obtained a positive answer to our
approximation question in a real-analytic situation for a new
class of functions $g$. By using a biholomorphic change of
coordinates, it is possible to assume that the first disk is the
standard disk $\{(z, \bar z)\, : \, z\in D\}$ and then apply an
approximation result of Nguyen, \cite{[2]}. In the present paper the same idea
of applying a biholomorphic map near the origin together with
already developed techniques in this area is used. We obtain
several new results of the form $[z^2, g^2 ; D]=C(D)$, one of them
being a generalization of the main result of \cite{[3]}, for new and larger classes of
functions $g$ (theorem 2.5).

\bigskip
{\bf Acknowledgement}.  The authors thank  Paul Beneker for a
stimulating discussion.

\section{An approximation result}
We agree on the following convention: all functions defined in a
neighborhood of the origin are of class $C^1$, even if we do not
mention this explicitely.
\begin{definition}
 Let  $g(z)$ be an {\it even}  function
 defined near
the origin with $g(z)=o(z)$.
 Suppose that there exists a polynomial $p(\zeta_1,
\zeta_2)$ such that for all functions $R(z)$  with $
R(z)=o(g(z))$
 both $$\im p(z, \bar
z+g(z)+R(z))>0$$
 and $$\im p(z, \bar z-g(z)+R(z))<0$$ 
hold for all
$z\ne 0$ sufficiently close to 0. \\ Then we say that $g$
satisfies the {\it polynomial condition} (with respect to $p$).
\end{definition}
\begin{examples} $\phantom{a}  $
\\
$\bullet \; \; \;$ If $m>1$, then for $g(z)=i|z|^m$ one can take
$p(\zeta_1, \zeta_2)=\zeta_1+\zeta_2$.
\\
$\bullet \; \;$ For the function $g(z)=a|z|^2+b\bar{z}^2$ with
$|b|<|a|$ one can take $p(\zeta_1, \zeta_2)=-ia\zeta_1+i\bar a
\zeta_2$. From this fact a version of the main result of [5]
follows.
\\
$\bullet \;\;$ The function $g(z)=|z|^2+\bar{z}
^2$ does not
satisfy the polynomial condition because it has non-zero zeroes.
\\
$\bullet \;\;$ The function $g(z)=z^3\bar z$  satisfies the
polynomial condition with respect to $p(\zeta_1,
\zeta_2)=-i\zeta_1^3+i\zeta_2^3$.
\end{examples}
\begin{lemma}
If $g$ satisfies the polynomial condition with respect to a
polynomial  $p$, then $g$ satisfies the  polynomial condition with
respect to the odd part of the polynomial  $p$.
\end{lemma}
\begin{proof}
Fix $R(z)=o(g(z))$  for the moment, then for $z\ne 0$ close to 0,
we have:
\begin{align} \im p(z,\bar{z}
+g(z)+R(z))&>0, \tag{a}\\
\im p(z, \bar{z}-g(z)-R(-z))&<0. \tag{b}
\end{align}
Replace  $z$ by $-z$ in (b) and use the fact that $g$ is even,
then also:
\begin{equation}
\im p(-z, -\bar{z}-g(z)-R(z))<0. \tag{c}
\end{equation}
 Now write $p$ as a sum of homogeneous
analytic polynomials, in other words $p=p_s+\cdots +p_n$ where
$p_j$ is homogeneous of degree $j$. Rewrite (c), for small $z\ne
0$, as:
$$\sum_{j=s}^n(-1)^{j+1}p_j(z, \bar z+g(z)+R(z))>0.$$
Combination with (a) shows that all terms with $j$ even in (a)
drop out.
\\ In a similar way these terms can be removed in the second part
of the polynomial condition.
\end{proof}
We need the following lemma which is without doubt well-known.
\begin{aux}
Let $F(w_1, w_2) $ be holomorphic near the origin, let $l\geq 2$
be an integer and let $F(w_1, w_2)=O(\| (w_1, w_2)\|^l)$.
\\
Let $A(w_1, w_2)$ be defined near the origin with $$A(w_1,
w_2)=O(\| (w_1, w_2)\|).$$
 Then sufficiently close to the origin
$$F(w_1, w_2+A(w_1, w_2))=F(w_1, w_2) + A(w_1, w_2)B(w_1, w_2),$$
with $B(w_1, w_2)=O(\| (w_1, w_2)\|^{l-1})$.
\end{aux}
\begin{proof}
As $F(w_1, w_2)$ is holomorphic near the origin,  $$H(w_1, w_2,
w_3)\, =
\begin{cases} \frac{F(w_1, w_3)-F(w_1, w_2)}{w_3-w_2}, &
\textnormal{if $w_3\ne w_2$,}
\\ \frac{\partial F}{\partial \zeta_2}(w_1, w_2), & \textnormal{if
$w_3=w_2$,}
\end{cases}$$
is holomorphic near the origin,  $H(w_1, w_2, w_3)=O(\|(w_1, w_2,
w_3)\|^{l-1})$ and
$$F(w_1, w_2+z)=F(w_1, w_2)+zH(w_1, w_2, w_2+z).$$
Since $A(w_1, w_2)=O(\| (w_1, w_2)\|)$ it follows that
$$F(w_1, w_2+A(w_1, w_2))=F(w_1, w_2) + A(w_1, w_2)B(w_1, w_2),$$
and  $$B(w_1, w_2)=H(w_1, w_2, w_2+A(w_1, w_2))=O(\| (w_1,
w_2)\|^{l-1}).$$
\end{proof}
\begin{theorem}
$\phantom{a}$
\\ $\bullet \; $ Let $F(w_1, w_2) $ be an odd holomorphic function
near the origin satisfying
 $F(w_1, w_2)=O(\| (w_1, w_2)\|^3)$ and let $f(z)=F(z, \bar
 z)$.
\\
$\bullet \; $ Suppose that $g$ satisfies the polynomial condition.
\\
$\bullet \; $ Let $h$  be defined near the origin with
$h(z)=o(g(z))$.
\\
Then for all
 disks $D$ about 0 with sufficiently small radius
$$[\; z^2, \; (\bar z+f(z)+g(z)+h(z))^2\; : \; D\; ]=C(D).$$
\end{theorem}
\begin{proof}
Let $X=\{\; (z^2, \; (\bar z+f(z)+g(z)+h(z))^2) \; : \; z\in D\;
\}.$
\\
The inverse image of $X$ under the map $\Pi : \Bbb{C}^2\to
\Bbb{C}^2$, defined by $\Pi(\zeta_1, \zeta_2)=(\zeta_1^2,
\zeta_2^2)$ consists of
\begin{align*}
 D_1&=\{\; (z, \bar z+f(z)+g(z)+h(z))\; :\; z\in D\;
\},
\\
D_2&=\{\; (-z, -(\bar z+f(z)+g(z)+h(z)))\; :\; z\in
D\; \} \\
 &=\{\; (z, \bar z+f(z)-g(z)-h(-z))\;
:\; z\in D\; \},
\\
 D_3&=\{\; (-z, \bar z+f(z)+g(z)+h(z))\; :\; z\in
D\; \},
\\
 D_4&=\{\; (z, -(\bar z+f(z)+g(z)+h(z)))\; :\; z\in
D\; \}
\\
&=\{\; (-z, \bar z+f(z)-g(z)-h(-z))\;
:\; z\in D\; \}.
\end{align*}
 Note that the condition on the existence of the polynomial $p$
implies that $g$ has no non-zero zeroes and that the two functions
$z^2$ and $(\bar z+f(z)+g(z)+h(z))^2$ separate the points of $D$
(if $D$ is sufficiently small).
\\
The techniques developed in the papers \cite{[8],[5]} on approximation on
disks give us:
\\
$\phantom{aaaaa} \phantom{\Longleftrightarrow} [\; z^2, \; (\bar
z+f(z)+g(z)+h(z))^2\; : \; D\; ]=C(D)$
\\ $\phantom{aaaa} \Longleftrightarrow  P(X)=C(X)$
\\ $\phantom{aaaa} \Longleftrightarrow  X$ is polynomially convex
\\ $\phantom{aaaa}
\Longleftrightarrow  D_1\cup D_2\cup D_3\cup D_4$ is polynomially
convex
\\ $\phantom{aaaa} \Longleftrightarrow  D_1\cup D_2$ is polynomially
convex.
\\
\\
We comment on these equivalences.
The first equivalence  is trivial.  Since $X$ is totally real except at the
origin, the second one follows from a theorem of O'Farrell, Preskenis and Walsh, \cite{[4]}. The
next equivalence is a consequence of a theorem of Sibony, \cite{[11]}, and the
last one is an application of Kallin's lemma using the polynomial
$p(\zeta_1, \zeta_2)=\zeta_1\cdot \zeta_2$.
\\
Later on we will also  use the following theorem of Wermer, 
\cite{[12]}. {\em If the function $F$ is of class $C^1$ near the origin in the
complex plane, with $F_{\bar z}(0)\ne 0$, then $[z,\, F : D]=C(D)$
if $D$ is a sufficiently small disk around 0.} This implies that
all disks $D_i$ are polynomially convex.
\\ For precise statements and use of these theorems, see \cite{[8]}, in particular the proof
of theorem 1.
\\
\\
Now let us show that $D_1\cup D_2$ is polynomially convex.
Consider the map  $G(w_1, w_2)=(w_1, w_2+F(w_1, w_2))$.  Since
$F(w_1, w_2)=O(\| (w_1, w_2)\|^3)$ it follows that $G$ is
biholomorphic near the origin (with inverse called $H$).
\\ Now $E_1=H(D_1)$ consists of points of the form $(z, q(z))$
where $q$ is of class $C^1$ near 0 and $q(0)=0$. Then there are $a$
and $b$ such that  $q(z)=az+b\bar z+r(z)$, where  $r(z)=o(z)$.
Applying $G$ we see  
\begin{equation}
(z, q(z)+F(z, q(z)))=(z, \bar
z+f(z)+g(z)+h(z)).\tag{$*$}
\end{equation} 
Since
$f(z)+g(z)+h(z)=O(z^3)+o(z)+o(z)$ and moreover $F(z, q(z))=O(z^3)$
 we infer that $q(z)=\bar z+r(z)$.  So $(*)$ translates into  
$$(z, \bar
z+r(z)+F(z, \bar z+r(z))) = (z, \bar z+f(z)+g(z)+h(z)).$$
 Applying the auxiliary lemma
  to this expression with $w_1=z, w_2=\bar z$ and 
$A(w_1, w_2)= r(w_1)$ we obtain: $$(z, \bar
z+r(z)+f(z)+r(z)B(z, \bar z))=(z, \bar z+f(z)+g(z)+h(z)).$$ 
It follows that $$r(z)=\frac{g(z)+h(z)}{1+B( z, \bar
z)}=g(z)+\frac{h(z)-g(z)B(z, \bar z)}{1+B(z, \bar z)}.$$
We conclude that
$E_1=H(D_1)$ consists of points $(z, \bar z+g(z)+R_1(z))$ in which
$R_1(z)=o(g(z))$ and is of class $C^1$. This last fact follows
from the definition of $B(w_1, w_2)$ in the proof of the auxiliary
lemma.
\\
Now $E_1$ is polynomially convex if $D$ is sufficiently small
(Wermer). \\
Similarly  $E_2=H(D_2)$ consists of points $(z, \bar
z-g(z)+R_2(z))$ in which $R_2(z)=o(g(z))$ and is of class $C^1$.
Also $E_2$ is polynomially convex if $D$ is sufficiently small.
Since $g$ satisfies the polynomial condition,
Kallin's lemma can be applied, showing that $E_1\cup E_2$ is polynomially convex. Applying
$G$ it follows that $D_1\cup D_2$ is polynomially convex for sufficiently small $D$.
\end{proof}
\begin{remark}
If $F(w_1, w_2)=f(w_1)=O(w_1^3)$  no computation is necessary
since
 the map
 $$G(w_1, w_2)=(w_1,
w_2+f(w_1))$$ has inverse $H(z_1, z_2)=(z_1, z_2-f(z_1))$ near the
origin. We now obtain directly $H(z, \bar z+f(z)+g(z)+h(z))=(z,
\bar z+g(z)+h(z))$ and similarly $H(z, \bar z+f(z)-g(z)-h(-z))=(z,
\bar z-g(z)-h(-z))$. Now use that $g$ satisfies  the polynomial
condition and proceed as before.
\end{remark}

\section{The polynomial condition for homogeneous functions}
Let $g$ satisfy the polynomial condition, then there is an odd
polynomial $p$ such that  
\begin{align}\im p(z, \bar
z+g(z)+R(z))>0\tag{1}\\ \intertext{and} \im p(z, \bar
z-g(z)+R(z))<0\tag{2}
\end{align}
 hold for all $z\ne 0$ sufficiently
close to 0 if $R(z) =o(g(z))$. As before $g$ is even, but instead
of $g(z)=o(z)$ we impose a stronger condition on this function:
$$g\; \hbox{is \textit{homogeneous} of degree $m>1$, i.e.}
$$ $$ g(tz)=t^mg(z)\; \hbox{ for}\; t>0$$
 (so in fact  $g$ is defined everywhere). Now write $p$ as a
sum of homogeneous analytic polynomials, $p=p_{2s-1}+\cdots +p_{2n-1}$ where all $p_k$ are
homogeneous of odd degree $k$. We assume first that $m$ is not an
odd integer. Let $n_0\leq n$ be maximal such that $2n_0-1<2s-2+m$.
\\  Taking for $R$ the zero function we
obtain:
\begin{multline*}
p(z,\bar z+g(z)) =p_{2s-1}(z, \bar z)+\cdots +p_{2n_0-1}(z, \bar
z)\\ + \frac{\partial p_{2s-1}}{\partial
\zeta_2}(z, \bar z)\cdot  g(z)+ O(|z|^{\alpha}),
\end{multline*} 
for some $\alpha>2s-2+m$. 
Now we restrict $z$ to the unit circle $\Gamma$, and obtain for $t>0$:
\begin{multline*}p(tz,
t\bar z+g(tz))  =t^{2s-1}p_{2s-1}(z, \bar
z)+\cdots +t^{2n_0-1}p_{2n_0-1}(z, \bar z)
\\ + t^{2s-2+m}\; \frac{\partial p_s}{\partial \zeta_2}(z, \bar z)\cdot g(z)+
O(t^{\alpha}).
\end{multline*}

Now take the imaginary part, divide by $t^{2s-1}$ and let $t$ tend
to 0. We obtain $\im p_{2s-1}(z, \bar z)\geq 0$. Similarly,
using the second condition on $g$, we obtain $\im p_{2s-1}(z, \bar
z)\leq 0$, hence $\im p_{2s-1}(z, \bar z)= 0$ for all $z\in \Gamma$
(hence for all $z\in \Bbb{C}$). Writing $p_{2s-1}(\zeta_1,
\zeta_2)=\sum_{k=0}^{2s-1}a_k\zeta_1^k\zeta_2^{2s-1-k}$ this means
that $a_k=\overline{a_{2s-1-k}}$ for all $k=0, \dots , 2s-1$. We
call such a polynomial {\it complex-symmetric}.
\\
Repeating this reasoning  we successively obtain:
$$ \im p_{2s+1}(z, \bar z)= 0, \dots , \im p_{2n_0-1}(z, \bar
z)= 0$$ and
\begin{equation}\im \frac{\partial p_{2s-1}}{\partial
\zeta_2}(z, \bar z)\cdot g(z) \geq 0. \tag{$*$}
\end{equation} Also in the case
that $m$ is an odd integer (1) and (2) in a similar way as above
lead to $(*)$.
 \\
 \\
Now suppose that for all $z\in \Gamma$ the inequality $(*)$ is
strict then we will show that the polynomial condition is satisfied for $g$ with respect to the polynomial
 $p_{2s-1}$. Indeed, if $R(z) =o(g(z))$ it follows for small $z\ne 0$:
\begin{multline*}p_{2s-1}(z,
\bar z+g(z)+R(z))
\\ =p_{2s-1}(z, \bar z) + \frac{\partial
p_{2s-1}}{\partial \zeta_2}(z, \bar z)\cdot  g(z)\cdot \bigl(
1+\frac{R(z)}{g(z)}\bigr) + O(|z|^{2s-3+2m}).
\end{multline*}
So for $z\in \Gamma$ and small $t>0$ it follows that:
\begin{multline*}\im p_{2s-1}(tz,t\bar z+g(tz)+R(tz))\\= \im
t^{2s-2+m}\Bigl( \frac{\partial p_{2s-1}}{\partial \zeta_2}(z,
\bar z)\cdot g(z)\cdot \bigl(
1+\frac{R(tz)}{g(tz)}\bigr)+O(t^{m-1})\Bigr).
\end{multline*}

Since  $\frac{R(tz)}{g(tz)}$ is uniformly small on $\Gamma$ if
$t>0$ is sufficiently small, the above expression is positive on
$\Gamma$ for small $t>0$. In other words:
$\im p_{2s-1}(z, \bar z+g(z)+R(z))>0$ if $z\ne 0$ is
sufficiently small. 
Also $\im p_{2s-1}(z, \bar z-g(z)+R(z))<0$ for small $z\ne
0$. So $g$ satisfies the polynomial condition with respect to
$p_{2s-1}$ and we proved:
\begin{theorem}
If $g$ is even and of class $C^1$ near the origin in the complex
plane, is homogeneous of order $m>1$ and satisfies $$\hbox{Im}\;
\frac{\partial p_{2s-1}}{\partial \zeta_2}(z, \bar z)\cdot g(z)>
0\;\; \hbox{ for all}\;\; z\in \Gamma ,$$ where $p_{2s-1}$ is a
homogeneous complex-symmetric polynomial of degree $2s-1$, then
$g$ satisfies the polynomial condition with respect to $p_{2s-1}$.
\end{theorem}
\begin{example}
An example of such a function is $g(z)=i\,
\overline{\frac{\partial p_{2s-1}}{\partial \zeta_2}(z, \bar z)}$,
where $p_{2s-1}$ is any homogeneous complex-symmetric polynomial of
degree $2s-1\geq 3$ ($s=1$ excluded because $g$ has to be
homogeneous of degree $m>1$) and such that $\frac{\partial
p_{2s-1}}{\partial \zeta_2}(z, \bar z)$ has no non-zero zeroes.
\end{example}
\begin{theorem} Let  
$g(z)=\sum_{k=-\infty}^{\infty} a_k\bar z^kz^{2m-k}
 $ with $m$ a positive integer. Suppose that 
 $ \sum_{k=-\infty}^{\infty}
|ka_k|<\infty$ and that one of the following increasingly weaker conditions is met:
$$\text{$\exists\ l\leq m$ such
that}\ |a_l|>\sum_{n\ne l}|a_n|, $$ 
 or 
$$\text{ $\exists\ l\leq m$ such
that}\ \sum_{n=1}^{\infty} \bigl| \frac{a_{l+n}}{a_l}+\frac{ \bar
a_{l-n}}{\bar a_l} \bigr| <1,$$
or 
$$\text{$\exists\ l\leq m$ such
that}\ \re \Bigl( 1+\sum_{n=1}^{\infty} \bigl(
\frac{a_{l+n}}{a_l}+\frac{ \bar a_{l-n}}{\bar a_l} \bigr)w^n
\Bigr)>0 \; \; \text{on} \;\;  |w|=1.$$

Then $g$ is an even  homogeneous $C^1$ function of degree $2m$
that satisfies the polynomial condition.
\end{theorem}
\begin{proof}
 Let $p(\zeta_1, \zeta_2)=\bar \alpha \zeta_1^{2m-2l+1}+\alpha
\zeta_2^{2m-2l+1}$ with $\alpha$ to be determined later (and with
$l\leq m$). Then for $z\in \Gamma$:
\begin{equation*}\begin{split}
&\frac{1}{2m-2l+1}\im\frac{\partial p}{\partial
\zeta_2}(z, \bar z)\cdot g(z)
=\im \sum_{k=-\infty}^{\infty} \alpha a_k\bar
z^{2m-2l+k}z^{2m-k}
\\
&=\im \{\sum_{k=-\infty}^{l-1}
\alpha a_k\bar z^{2m-2l+k}z^{2m-k} 
  + \alpha a_l |z|^{4m-2l}  +
\sum_{k=l+1}^{\infty} \alpha a_k\bar z^{2m-2l+k}z^{2m-k} \}
\\
&=\im\{
  \alpha a_l |z|^{4m-2l }  + \sum_{n=1}^{\infty}
\bigl( \alpha a_{l+n}-\bar \alpha \bar a_{l-n} \bigr) \bar
z^{2m-l+n}z^{2m-l-n }\}
\\
&= \im\{ i|a_l||z|^{4m-2l}\Bigl( 1+
\sum_{n=1}^{\infty} \bigl(\frac{a_{l+n}}{a_l}+\frac{ \bar
a_{l-n}}{\bar a_l} \bigr) \bigl( \frac{\bar
z}{z}\bigr)^n\Bigr)\}.\end{split}
\end{equation*}
In the last equality we chose $\alpha =i\frac{|a_l|}{a_l}$. 
The final expression has positive imaginary part if the third condition in the statement of the theorem is satisfied.
\end{proof}
\begin{remarks}
  This result includes the more restricted
case of polynomials
 $g(z)=\sum_{k=0}^{2m} a_k\bar{z}^kz^{2m-k}$ in $z$ and $\bar z$,
 for which there exists $0\leq l\leq m$ such
that $|a_l|>\sum_{k\ne l}|a_k|$, essentially studied by Nguyen,
\cite{[2]}, and applied
  in a real-analytic setting by Nguyen and De Paepe, \cite{[3]}. The
  condition on the coefficients here is more general. For instance
  if $m=1$ the condition is valid if $\bigl|
  \frac{a_2}{a_1}+\frac{\bar{a}_0}{\bar{a}_1}\bigr|<1$, which is
  certainly the case for (but is not equivalent to) $|a_1|>|a_0|+|a_2|$.
  \end{remarks}
\begin{example}
 Applying theorem 3.3 and theorem 2.5 we obtain a result from \cite{[9]}:
 $$[\, z^2, \, \bar z^2 + z^3 ; \, D\, ]=[\, z^2, \, (\bar z+\frac{1}{2}\frac{z^3}{\bar z} + \hbox{h.o.t.} )^2\, ; \,
D\, ]=C(D).$$
\end{example}

\section{Another use of a biholomorphic map}
In theorem 2.5 it was fruitful to apply a biholomorphic map in
order to show polynomial convexity. This idea can be used in other
situations as well. For instance, suppose that $g$ is of class $C^1$
near 0, $g(0)=0, g_z(0)=0, g_{\bar z}(0)=1$ and such that $z^2$
and $g^2$ separate points near 0.
\\
Also suppose $F$ is defined near the origin,  holomorphic, and odd,
with $F(w_1, w_2)=O(\| (w_1, w_2)\|^3)$. Then $z^2$ and $(g+F(z,
g))^2$ separate points near 0 and $[\, z^2, \, (g+F(z, g))^2\, ;
\, D\, ]\subset [\, z^2, \, g^2\, ; \, D\, ]$. So  $[\, z^2, \,
g^2;  D\, ]\ne C(D)$ implies $[\, z^2, \, (g+F(z, g))^2 ; D\, ]\ne
C(D)$. This is the contents of the proof of theorem 2 in \cite{[9]}. But
more is true.
\begin{theorem}
With notation as above and for sufficiently small $D$:
$$[\, z^2, \, g^2;  D\, ]=C(D)\Longleftrightarrow [\, z^2, \, (g+F(z,
g))^2 ; D\, ]=C(D).$$
\end{theorem}
\begin{proof}
Let $X=\{\, (z^2, \, g(z)^2) \, : \, z\in D\, \}$, furthermore if
we  let \\
$Y=\{\, (z^2, \, (g(z)+F(z, g(z)))^2) \, : \, z\in D\, \}$, then,
using the biholomorphic map $G(w_1, w_2)=(w_1, w_2+F(w_1, w_2))$
in the fourth equivalence, we obtain for sufficiently small $D$:
\\
$\phantom{\Longleftrightarrow}\; [\, z^2, \, g^2; D\, ]=C(D)$
\\
$ \Longleftrightarrow P(X)=C(X)$
\\
$ \Longleftrightarrow X $ is polynomially convex
\\
$ \Longleftrightarrow \{\, (z, g(z))\, :\, z\in D\, \}\cup  \{\,
(z, -g(-z)) \, :\, z\in D\, \}$ is pcx
\\
$ \Longleftrightarrow \{\, (z, g(z)+F(z, g(z)))\, :\, z\in D\,
\}$\\$\phantom{aaaaaaaaaa} \cup \{\, (z, -g(-z)+F(z, -g(-z))) \,
:\, z\in D\, \}$ is pcx
\\
$ \Longleftrightarrow Y $ is polynomially convex
\\
$ \Longleftrightarrow P(Y)=C(Y)$
\\
$\Longleftrightarrow [\, z^2, \, (g+F(z, g))^2 : D\, ]=C(D)$.
\end{proof}
\begin{question}
Is $[\, z^2, \, (g+F(z, g))^2\, : \, D\, ]=[\, z^2, \, g^2\, ; \,
D\, ]$ for all $g$ and $D$ as above?
\end{question}
\begin{example}
In the case $F(w_1, w_2)=f(w_1)=O( w_1^3)$ the answer to the
question is yes:
$$[\, z^2, \, (g+f)^2\, : \, D\, ]=[\, z^2, \, g^2\, : \,
D\, ].$$
 Indeed, since $zf, f^2$ and $\frac{f}{z}$ are even analytic
functions, they belong to $A=[\, z^2, \, (g+f)^2\, : \, D\, ]$.
Also $z^2(g+f)^2\in A$, thus (since the real part of $z^2(g+f)^2$ is
non-negative near the origin)  $z(g+f)\in A$, hence $zg\in A$. Also
$(g+f)^2=g^2+2(zg)\cdot \frac{f}{z}+f^2\in A$, therefore $g^2\in A$. Hence
$A=[\, z^2, \, g^2\, : \, D\, ].$
\end{example}
\begin{example}
 A second situation where the answer is yes occurs when $F(w_1,
w_2)$ has the form $w_2G(w_1^2, w_2^2)$ where $G$ is holomorphic
near the origin with $G(0, 0)=0$. Then  $(g+F(z, g))^2$ can be
written as $g^2+g^2H(z^2, g^2)$ with $H(0, 0)=0$. The map
$$(w_1, w_2)\mapsto (w_1, w_2 +w_2H(w_1, w_2))$$
is biholomorphic near the origin and maps the pair $(z^2, g^2)$ to
$(z^2, (g+F(z, g))^2)$. This shows that the algebra generated by
$z^2$ and $g^2$ on  a small $D$ equals the algebra generated by $z^2$ and
$(g+F(z, g))^2$  on $D$.
\end{example}
\bigskip
\section{Appendix}
In this appendix we keep the setting of section 3 and see what
can be said when $\im\frac{\partial
p_{2s-1}}{\partial \zeta_2}(z, \bar z)\cdot g(z)$ has zeroes on
$\Gamma$. Under stronger conditions on $F$, $g$, and $h$,
we obtain the following approximation result.
\\
\begin{theorem} Let $F(w_1, w_2) $ be an odd holomorphic function near the origin satisfying
 $F(w_1, w_2)=O(\| (w_1, w_2)\|^5)$ and set $f(z)=F(z, \bar
 z)$. Suppose that the following conditions are met:
\begin{itemize}
\item The function $g\in C^1$ is even and homogeneous
of degree $m>3$.
\item There are homogeneous complex-symmetric
polynomials $p_{2s-1}$ and $p_{2s+1}$ of degree $2s-1$,
respectively $2s+1$, such that  $$\im \frac{\partial
p_{2s-1}}{\partial \zeta_2}(z, \bar z)\cdot g(z)\geq 0,\;\; \hbox{
for all}\;\; z\in \Gamma, $$ and
$$ \im\frac{\partial
p_{2s+1}}{\partial \zeta_2}(z, \bar z)\cdot g(z)> 0$$ for
all $ z\in \Gamma$ where $\im \frac{\partial
p_{2s-1}}{\partial \zeta_2}(z, \bar z)\cdot g(z)= 0$.
\item  $h$ is defined near the origin and
$h(z)=o(z^2g(z))$.
\end{itemize}
Then for all
 disks $D$ centered at 0 with sufficiently small radius
$$[\; z^2, \; (\bar z+f(z)+g(z)+h(z))^2\; ; \; D\; ]=C(D).$$
\end{theorem}
\begin{proof}
We will follow the line of the proofs of the auxiliary lemma 2.4 and of theorem 2.5, as well as the notation. We
see that $B(z, \bar z)=O(z^4)$ since $F(w_1, w_2)=O(\| (w_1,
w_2)\|^5)$. From this fact and $h(z)=o(z^2g(z))$ it follows that
$R_1(z), R_2(z) = o(z^2g(z))$.
\\
\\
Let $N$ be the set of points $z\in \Gamma$ where $$\im
\frac{\partial p_{2s-1}}{\partial \zeta_2}(z, \bar z)\cdot g(z)=0.
$$
 Now assume $$\im
\frac{\partial p_{2s+1}}{\partial \zeta_2}(z, \bar z)\cdot
g(z)>0\; \;  \hbox{for all}\; \;  z\in N.$$
 Then there is ${\lambda}_0 >0$ and 
$\delta >0$ such that for all $z\in \Gamma$ and $0<\lambda\leq
{\lambda}_0$:
$$\im\frac{\partial \bigl( p_{2s-1}+\lambda p_{2s+1}\bigr) }{\partial \zeta_2}(z, \bar z)\cdot
g(z) \ge \lambda\delta. $$ 
Indeed, for $z\in \Gamma$, let
$$f_0(z)=\im \frac{\partial p_{2s-1}}{\partial \zeta_2}(z, \bar
z)\cdot g(z),\qquad
 f_1(z)=\im\frac{\partial
p_{2s+1}}{\partial \zeta_2}(z, \bar z)\cdot g(z).$$
\\ Let $0<2\delta = \hbox{inf}_{\; z\in N \;}  f_1(z)$,
$U$  a neighbourhood of $N$ in $\Gamma$ such that $
\hbox{inf}_{\; z\in U \;}  f_1(z)\geq \delta$ and $\epsilon =
\hbox{inf}_{\; z\in \Gamma \setminus U \;}  f_0(z)>0$. \\
If we take $0< \lambda \leq \lambda_0=\hbox{min}\{
\frac{\epsilon/2}{\|f_1\|_{\Gamma}}\; , \;
\frac{\epsilon}{2\delta}\}$, then $f_0+\lambda f_1\geq \lambda
\delta $ on $\Gamma$.

\bigskip
Now for $m>3$ and $R(z)=o(z^2g(z))$ we have:
\begin{multline*}(p_{2s-1}+p_{2s+1})(z, \bar z+g(z)+R(z))
= p_{2s-1}(z, \bar z)+p_{2s+1}(z, \bar z)\\
 +\bigl( \frac{\partial p_{2s-1}}{\partial \zeta_2}(z, \bar
z) + \frac{\partial p_{2s+1}}{\partial \zeta_2}(z, \bar z)\bigr)
\cdot g(z)\cdot \bigl( 1+\frac{R(z)}{g(z)}\bigr)  +O(|z|^{2s-3+2m}).
\end{multline*}
\\
So for  $z\in
\Gamma$ and small $t>0$  one has \begin{multline*}
 \im (p_{2s-1}+p_{2s+1})(tz, t\bar
z+g(tz)+R(tz))
\\
 = t^{2s-2+m}\hbox{Im}\;\bigl[ \bigl( \frac{\partial
p_{2s-1}}{\partial \zeta_2}(z, \bar z)+ t^2\frac{\partial
p_{2s+1}}{\partial \zeta_2}(z, \bar z)\bigr)\\
 \cdot g(z)\cdot
\bigl( 1+t^2\cdot \frac{z^2R(tz)}{(tz)^2g(tz))}\bigr)+ O(t^{m-1})\bigr]
  \geq \frac{1}{2}t^{2s+m}\delta, 
\end{multline*} since  $\frac{z^2R(tz)}{(tz)^2g(tz)}$ is uniformly small on $\Gamma$ if $t>0$ is sufficiently small.
\\ So $\im
(p_{2s-1}+p_{2s+1})(z, \bar z+g(z)+R(z))>0$,  and similarly  $\im (p_{2s-1}+p_{2s+1})(z, \bar z-g(z)+R(z))<0$ if $z$ is sufficiently small.
 Now proceed as in the proof of theorem 2.5.
\end{proof}
\bigskip
\bibliographystyle{amsplain}

\end{document}